\documentclass{amsart}
\newtheorem{Theorem}{Theorem}[section]
\newtheorem{Definition}[Theorem]{Definition}
\newtheorem{Proposition}[Theorem]{Proposition}

\newtheorem{Lemma}[Theorem]{Lemma}
\newtheorem{Corollary}[Theorem]{Corollary}
\theoremstyle{remark}

\newtheorem{Example}[Theorem]{Example}

\newcommand\eps{\varepsilon}

\newcommand\ovr{\overline}

\newcommand\Om{\Omega}
\newcommand\al{\alpha}

\newcommand\dl{\delta}

\newcommand\bd{\partial}
\newcommand\lm{\lambda}

\newcommand\si{\sigma}

\newcommand\sm{\setminus}
\newcommand\sbs{\subset}

\newcommand\wtl{\widetilde}

\newcommand\re{{\mathbf {Re\,}}}

\newcommand\be{\begin{enumerate}}
\newcommand\ee{\end{enumerate}}
\newcommand\bT{\begin{Theorem}}
\newcommand\eT{\end{Theorem}}
\newcommand\bP{\begin{Proposition}}
\newcommand\eP{\end{Proposition}}
\newcommand\bD{\begin{Definition}}
\newcommand\eD{\end{Definition}}
\newcommand\bE{\begin{Example}}
\newcommand\eE{\end{Example}}
\newcommand\bL{\begin{Lemma}}
\newcommand\eL{\end{Lemma}}
\newcommand\bC{\begin{Corollary}}
\newcommand\eC{\end{Corollary}}
\newcommand\A{{\mathcal A}}

\newcommand\oD{\ovr{\mathbb D}}
\newcommand\aD{\mathbb D}

\newcommand\aR{\mathbb R}

\newcommand\aC{\mathbb C}
\newcommand\aCP{\mathbb {CP}}

\begin{document}
\title{Pluricomplex Green functions on manifolds}
\author{Evgeny A. Poletsky}
\begin{abstract}  In this paper we prove the basic facts for pluricomplex Green functions on manifolds. The main goal is to establish properties of complex manifolds that make them analogous to relatively compact or hyperconvex domains in Stein manifolds.
\thanks{The author was partially supported by a grant from Simons Foundation.}
\end{abstract} 
\address{ Department of Mathematics,  215 Carnegie Hall,
Syracuse University,  Syracuse, NY 13244}
\maketitle
\section{Introduction}\label{S:1}
\par The pluricomplex Green functions are one of the most popular biholomorphic invariants in several complex variables. Their properties on relatively compact domains in  Stein manifolds or hyperconvex relatively compact domains in Stein manifolds got an excellent exposition in the books of Klimek \cite{K2} and Jarnicki and Pflug \cite{JP} combined with the papers \cite{D1} of Demailly and \cite{B} of B\l ocki. In these cases the attainable boundaries are present. The goal of this paper is to describe complex manifolds, where the pluricomplex Green functions have the same basic properties as on domains above but the boundaries are not visible.
\par Important basic features of relatively compact domains in Stein manifolds are:\be
\item Carath\'eodory and Kobayashi hyperbolicity;
\item separation of points by bounded holomorphic functions;
\item existence of the Bergman metric;
\item existence of the pluricomplex Green functions with the strict logarithmic pole;
\item existence of continuous bounded strictly plurisubharmonic functions.
\ee
\par In Section \ref{S:slp} we show that the existence of a continuous bounded strictly plurisubharmonic function on a complex manifold $M$ allows us to preserve all these properties but Carath\'eodory hyperbolicity. This condition already appeared in the literature. Chen and Zhang showed in \cite{CZ} that manifolds with this property possess the Bergman metric and Shcherbina and the author proved in \cite{PSh} that their points can be separated by bounded continuous plurisubharmonic functions.
\par Important basic features of hyperconvex relatively compact domains in Stein manifolds are:\be
\item continuity of the pluricomplex Green functions in both variables;
\item the pluricomplex Green functions are exhaustion functions;
\item completeness of the Bergman and Kobayashi metrics;
\item existence of a smooth strictly plurisubharmonic negative exhaustion function.
\ee
\par In \cite{C1} Chen introduced hyperconvex manifolds as manifolds with negative strictly plurisubharmonic exhaustion function. In Section \ref{S:hcm} we show that this definition is equivalent to requiring the existence of a negative plurisubharmonic exhaustion function and a continuous bounded strictly plurisubharmonic function and show that all the features listed above still hold.
\par This paper was conceived  when the author was visiting N. Shcherbina at the University of Wuppertal in 2018. Many of its ideas and proofs are based on fruitful discussions with N. Shcherbina and our joint work. I am grateful to the university and Nikolay for their hospitality and encouragement.
\par The author is also grateful to the referees whose comments and suggestions significantly improved the exposition.
\section{Definitions and basic properties}
\par Throughout this paper $\aD(0,r)=\{\zeta\in\aC:\,|\zeta|<r\}$, $\aD=\aD(0,1)$ and $B(w,r)=\{z\in\aC^m:\,\|z-w\|<r\}$.
\subsection{Definitions:} Let $M$ be a connected complex manifold. A function on $M$ has a {\it logarithmic pole} at $w\in M$ if $u(z)-\log\|z-w\|$ is bounded above on some coordinate neighborhood of $w$. When $M$ is a domain in $\aC^n$ Klimek in \cite{K1} and, independently under the name ``invariant Green function'', the author in \cite{PS},  introduced for $z,w\in M$ the {\it pluricomplex Green function}
\[g_M(z,w)=\sup u(z),\]
where the supremum is taken over all negative plurisubharmonic functions $u$ with logarithmic pole at $w$. (Here we assume that $u\equiv-\infty$ is a plurisubharmonic function.) This definition word to word can be used to define the pluricomplex Green function on complex manifolds. In \cite{D1} Demailly used it for relatively compact domains in Stein manifolds.
\par Let $\A_r(M)$ be the space of  {\it analytic disks} in $M$, i. e., continuous mappings $f$ of the closed disk $\ovr\aD(0,r)$ into $M$ holomorphic on $\aD(0,r)$. We will write $\A(M)$ for $\A_1(M)$. Given $w\in M$ consider on $\A(M)$ the functional
\[H_w(f)=\sum_{f(\zeta_j)=w}\log|\zeta_j|.\]
When $M$ is a domain in $\aC^n$ the author of this paper have shown in \cite{PS} that \[g_M(z,w)=\inf\{H_w(f):\,f\in\A(M), f(0)=z\}.\]
Later L\'arusson and Sigurdsson in \cite{LS} proved this formula for all complex manifolds.
\subsection{Basic properties:} First of all, if $M$ and $N$ are complex manifolds and $F:\,M\to N$ is a holomorphic mapping, then $g_M(z,w)\ge g_N(F(z),F(w))$. In particular, $g_M$ is monotonically decreasing in $M$ and is biholomorphically  invariant.
\par Secondly, if $M$ is a ball $B(w_0,r)$ in $\aC^n$, then the function $g_M(z,w_0)=\log(\|z-w_0\|/r)$. Consequently, if $U\sbs M$ is a coordinate neighborhood of $w\in M$, then there is a constant $c$ depending only on $U$ such that $u(z)-\log\|z-w\|\le c$ near $w$ when $u$ is a plurisubharmonic function on $M$ with a logarithmic pole at $w$. Thus the upper regularization of $g_M(\cdot,w)$ also has a logarithmic pole at $W$ and, since it is plurisubharmonic, it is equal to $g_M(z,w)$ and we see that $g_M(z,w)$ is plurisubharmonic in $z$.
\par The function $g_M(z,w)$ is also maximal in $z$ outside $w$ in the following sense: if $G\sbs M$ is a domain whose closure does not contain $w$ and $v$ is a negative plurisubharmonic function on a neighborhood $U$ of $\ovr G$ such that $v(z)\le g_M(z,w)$ on $\bd G$, then $v(z)\le g_M(z,w)$ on $G$. Indeed, take the function $v_1(z)$ equal to $g_M(z,w)$ on $M\sm G$ and to $\max\{g_M(z,w),v(z)\}$ on $G$. This function will be negative and plurisubharmonic on $M$ and $v_1(z)=g_M(z,w)$ near $w$. Thus $g_M(z,w)\ge v(z)$ on $G$.
\subsection{Connections to invariant distances:} Let us introduce the {\it Carath\'eodory function}
\[c_M(z,w)=\sup\{\log|f(z)|:\,f:\,M\to\aD\text{ is holomorphic and } f(w)=0\}\]
and the {\it Kobayashi function}
\[k_M(z,w)=\inf\{\log|\zeta|:\,f\in\A(M), f(0)=z, f(\zeta)=w\}.\]
Clearly, $c_M(z,w)\le g_M(z,w)\le k_M(z,w)$.
\par The Carath\'eodory and Kobayashi functions are symmetric, the former function is continuous, the latter one is upper semicontinuous and both have a logarithmic pole at $w$. Hence if the function $k_M(z,w)$ is plurisubharmonic in $z$, then it is equal to $g_M(z,w)$.
\par According to Lempert  \cite{L1} on smooth strongly convex domains the Kobayashi and Carath\'eodory distances coincide. Hence in this case (see \cite{PS} or \cite{K1}) $c_M(z,w)=g_M(z,w)=k_M(z,w)$.
\par Almost nothing is known about these equalities for other classes of manifolds. The only result in this direction known to the author is the identity $g_M(z,w)=k_M(z,w)$ proved by Krushkal in \cite{Kr} when $M$ is the Teichm\"uller space $T_{g,n}$.
\par The pluricomplex Green function found applications to the Bergman metric, in particular, to show its completeness. The references to the papers with such applications when $M$ is a domain in $\aC^n$ can be found in \cite{B}, \cite{C2} and \cite{JP}. We will recall the results for manifolds later.
\section{Pluricomplex Green function with strict logarithmic pole}\label{S:slp}
\subsection{Pluri-Greenian manifolds:} A domain in $\aR^n$ is called {\it Greenian} (see \cite{AG}) if it has the Green functions with the right pole at any point. For complex manifolds the right pole is a {\it strict logarithmic pole.} A function $u$ on $M$ has a strict logarithmic pole at $w\in M$ if there are a coordinate neighborhood $U$ of $w$ and constants $c_1$ and $c_2$ such that
\[\log\|z-w\|+c_1\le u(z)\le\log\|z-w\|+c_2\]
on $U$.
\par Following the tradition to add the word ``pluri'' to pluripotential analogs of notions in the potential theory we will call a complex manifold $M$ {\it pluri-Greenian} if $g_M(z,w)$ has a strict logarithmic pole at any $w\in M$. On pluri-Greenian manifolds the pluricomplex Green functions $g_M(z,w)$ are locally bounded outside $w$ and, therefore, are maximal in the sense that $(dd^c_zg_M(z,w))^m\equiv0$ on $M\sm\{w\}$, where $m=\dim M$. Hence $(dd^c_zg_M(z,w))^m=(2\pi)^n\dl_w$.
\par Let us give some examples of complex manifold that are not pluri-Greenian.
\subsection{Parabolic manifolds:} A complex manifold $M$ is {\it parabolic} if any bounded above plurisubharmonic function on $M$ is constant. On these manifolds $g_M(z,w)\equiv-\infty$. Typical examples are compact complex manifolds with a locally pluripolar set removed. For more examples see \cite{AS}.
\par An $S$-manifold, introduced in \cite{AS}, is a complex manifold $M$ with an unbounded plurisubharmonic exhaustion function that is maximal outside a compact set in $M$. Demailly proved in \cite{D} that any $S$-manifold is parabolic. The converse is unknown. An interesting problem is: if $M$ is $S$-parabolic does it have an unbounded plurisubharmonic exhaustion function $u$ such that $(dd^cu)^m$ is an atom?
\subsection{Non-empty cores:}  Let $PSH^{cb}(M)$ be the space of bounded above continuous plurisubharmonic functions on $M$. The  set $\mathbf{c}(M)$ of all
points $w\in M$, where every function of $PSH^{cb}(M)$ fails to be smooth
and strictly plurisubharmonic near $w$, is called the {\it core} of $M$. This notion was introduced by Harz--Shcherbina--Tomassini in \cite{HST1}.
\par If the core is non-empty, then by \cite[Theorem 4.8]{PSh} it can be decomposed into the disjoint union of closed sets $E_j$, $j\in J$, such that each of these sets has no isolated points and has the following Liouville property: every function $\phi\in PSH^{cb}(M)$ is constant on each of the sets $E_j$. Thus any negative continuous plurisubharmonic function cannot have a strict logarithmic pole at any $w\in\mathbf{c}(M)$. It is unknown what happens if the assumption of continuity is removed.
\par For example, if $M=\aD_\xi\times\mathbb C_\zeta$, then $\mathbf{c}(M)=M$
and $E_\xi=\{(\xi,\zeta):\,\xi\in \aD, \zeta\in\aC\}$. The pluricomplex Green function $g_M((\xi,\zeta),(\xi_0,\zeta_0))=g_{\aD}(\xi,\xi_0)$ does not have a strict logarithmic pole at any point of $M$.
\par The close connection between the notions of a core and a strict logarithmic pole was demonstrated in the following theorem proved in \cite[Theorem 3.2]{PSh}.
\bT\label{T:3.2} A point $w\not\in \mathbf{c}(M)$ if and only if  there
is a negative continuous plurisubharmonic function $v$ on $M$ with a strict logarithmic pole at $w$.\eT
\par It follows from this theorem that the function $g_M(z,w)$ has a strict logarithmic pole at any $w\not\in \mathbf{c}(M)$ and this allows us to give an example of a domain in $\aC^n$, $n\ge2$, that has a strict logarithmic pole at some point but is not pluri-Greenian. Recall that if $M$ is a domain in the complex plane and $g_M$ has a strict logarithmic pole at some point, then $M$ is Greenian.
\bE Take a Fatou--Biberbach domain $\Om$ in $\aC^2$ and a point $w_0\in\bd\Om$. Let $B$ be a ball centered at $w_0$ whose boundary contains a point $w_1$ in the complement of $\ovr\Om$ in $\aC^2$. Let $\Om'=\Om\cup B$. We may assume  that $w_0=(-1,0)$ and $w_1=0$. Let $\rho(z)=\|z-w_0\|^2+\re z_1-1$, where $z=(z_1,z_2)$. Let us take $\dl>0$ such that the set $E=B\cap\{z:\,-\dl<\re z_1<0\}$ does not meet $\Om$. Define $u(z)=\max\{\rho(z),-\dl\}$ on $E$ and $-\dl$ on $\Om'\sm E$. Note that $\rho(z)<-\dl$ on $E=B\cap\{z:\,\re z_1=-\dl\}$ and $\rho(z)>-\dl$ on some neighborhood $U$ of 0 in $B$. Hence $u$ is plurisubharmonic and any point $z\in U$ does not belong to $\mathbf{c}(M)$. So at any $w\in U$ the function $g_{\Om'}(z,w)$ has a strict logarithmic pole but if $w\in\Om$, then $g_M(z,w)\equiv-\infty$.
\eE
\par However, it was shown in \cite{GHH} that a strongly pseudoconvex domain in $\aC^m$ cannot contain Fatou--Biberbach domains.
\subsection{Locally uniformly pluri-Greenian complex manifolds:} We introduce {\it locally uniformly pluri-Greenian} complex manifolds $M$, where every point $w_0\in M$ has a coordinate neighborhood $U$ with the following property:  there is an open set $W\sbs U$ containing $w_0$ and a constant $c$ such that $g_M(z,w)\ge\log\|z-w\|+c$ on $U$ whenever $w\in W$.
\par If $M$ is a ball $B(w_0,r)$ in $\aC^m$, then $g_M(z,w_0)=\log(\|z-w_0\|/r)$. Since $g_M$ is monotonic in $M$, it follows that if $M$ is a bounded domain in $\aC^n$, then $g_M(z,w)\ge\log(\|z-w\|/r)$,
where $r$ is the radius of circumscribed ball of $M$ centered at $w$. Hence bounded domains in $\aC^n$ are locally uniformly pluri-Greenian.
\par The following important result, proved in \cite{D1} for relatively compact domains in Stein manifolds, stays valid for locally uniformly pluri-Greenian manifolds (see \cite[Lemma 6.1]{P}).
\bL\label{L:kl} If $M$ is a locally uniformly pluri-Greenian complex manifold  and $w_0\in M$, then for any $\eps>0$ and any neighborhood $X$ of $w_0$ there is a neighborhood $Y$ of $w_0$ such that
\[1-\eps\le\frac{g_M(z,w_0)}{g_M(z,w)}<1+\eps\]
whenever $w\in Y$ and $z\in M\sm X$.
\eL
\par As immediate corollary we get
\bC\label{C:ucg} If $M$ is a locally uniformly pluri-Greenian complex manifold, then the pluricomplex Green function $g_M(z,w)$ is continuous in $w$.
\eC
\par In \cite{P} a special norm was introduced on the space $PSH^-(M)$ of negative plurisubharmonic functions on $M$. If we assign to each point $w\in M$ the function $\wtl g_M(z,w)=\al^{-1}(y)g_M(z,w)$, where $\al(y)$ is the norm of $g_M(z,w)$, then we obtain an imbedding $\Phi$ of $M$ into $PSH^-(M)$.   By Corollary \ref{C:ucg} this imbedding is a homeomorphism of $M$ onto $\Phi(M)$ when $M$ is a locally uniformly pluri-Greenian complex manifold and the properties of the norm imply that the closure $\wtl M$ of $\Phi(M)$ is compact. Moreover, the following result holds (see \cite{P} for details):
\bT Let $M$ and $N$ be locally uniformly pluri-Greenian complex manifolds. Then any biholomorphic mapping $F:\,M\to N$ extends to a homeomorphism of $\wtl M$ onto $\wtl N$. \eT
\par We computed in \cite{P} the pluripotential compactification for a ball, smooth strongly convex domains and a bidisk.  In all cases the limits of $\wtl g_M(z,w_j)$ as $\Phi(w_j)\to\bd\wtl M$ are maximal on $M$ and are the scalar multiples of the pluriharmonic Poisson kernels computed in \cite{D1} and \cite{BPT}. In the two first cases the set $\bd\wtl M$ coincides with the Euclidean boundary while in the case of a bidisk it is the product of a circle and a 2-sphere.
\subsection{Complex manifolds with bounded strictly plurisubharmonic function:}
\par Following Narasimhan in \cite{Nar} (see also \cite[Definition 5.20]{D2}) we say that a plurisubharmonic function $u$ on a complex manifold $M$ is {\it strictly plurisubharmonic } if for every point $w\in M$ there is a coordinate neighborhood $U$ containing $w$ and a number $\eps>0$ such that the function $u(z)-\eps\|z-w\|^2$ is plurisubharmonic on $U$. When $u$ is $C^2$-smooth this definition coincides with the standard definition of strict plurisubharmonicity.
\par For transition from continuous strictly plurisubharmonic functions to smooth strictly plurisubharmonic functions we will need the following result of Richberg in \cite{R} (see also \cite[Theorem 5.21]{D2}).
\bT\label{T:rt} Let $u$ be a continuous strictly plurisubharmonic function on a complex manifold $M$. For any continuous positive function $\lm$ on $M$ there is a smooth strictly plurisubharmonic function $v$ on $M$ such that $u\le v\le u+\lm$.
\eT
\par The following theorem claims that for a complex manifold to be pluri-Greenian is almost equivalent to the possession of a continuous bounded strictly plurisubharmonic function.
\bT\label{T:gr} If a connected complex manifold $M$ has a bounded continuous strictly plurisubharmonic function, then it is pluri-Greenian. Conversely, if for any $w\in M$ there is a negative continuous plurisubharmonic function $v$ on $M$ with a strict logarithmic pole at $w$, then $M$ has a bounded continuous strictly plurisubharmonic function.
\eT
\begin{proof} If there is a bounded continuous strictly plurisubharmonic function $u$ on $M$, then by Theorem \ref{T:rt} there is a smooth bounded strictly plurisubharmonic function $\wtl u$ on $M$. Hence ${\mathbf c}(M)=\emptyset$ and by Theorem \ref{T:3.2} $g_M$ has a strictly logarithmic pole at all points and, therefore, $M$ is pluri-Greenian.
\par If for any $w\in M$ there is a negative continuous plurisubharmonic function $v$ on $M$ with a strict logarithmic pole at $w$, then by Theorem \ref{T:3.2} the core $\mathbf{c}(M)$ is empty and again for any $w\in M$ there is a continuous bounded plurisubharmonic function on $M$ that is strictly plurisubharmonic near $w$.
\par Let us cover $M$ by countable family of open sets $U_j$ such that for each $j$ there is a continuous bounded plurisubharmonic function $u_j$ on $M$ that is strictly plurisubharmonic on $U_j$. Let $\al_j$ be the sup-norm of $u_j$ on $M$ and  $u=\sum_j2^{-j}\al_j^{-1}u_j$. Then $u$ is a bounded strictly plurisubharmonic function on $M$.  \end{proof}

\par By the following  theorem manifolds with a bounded continuous strictly plurisubharmonic function  are locally uniformly pluri-Greenian. Its proof is a minor elaboration on the proof of Theorem 3.2 in \cite{PSh}.
\bT\label{T:bsplug} If a complex manifold $M$ has a bounded continuous strictly plurisubharmonic function, then it is locally uniformly pluri-Greenian.
\eT
\begin{proof} By Theorem \ref{T:rt} the existence of bounded continuous strictly plurisubharmonic function on $M$ implies the existence of smooth bounded strictly plurisubharmonic function $u$ on $M$. We may assume that $u<0$ on $M$.
\par Let $w_0\in M$. By Lemma 3.1 in \cite{PSh} in some coordinate neighborhood $V$ of $w_0$ there is a ball $B$ of radius $r$ centered at $w_0$ and a pluriharmonic function $h$ on  $V$ such that
$u-h>0$ near $\bd B$ and $a=u(w_0)-h(w_0)<0$. There is another ball
$U=B(w_0,s)\sbs\sbs B$, $s>0$, such that $u-h<a/2$ on $\ovr U$.
The pluricomplex Green function $g_B(z,w_0)=\log(|z-w_0|/r)$ on $B$.  We
take a constant $d>0$ such that $d\log(s/r)>a/4$. Then $dg_B(z,w_0)>a/4$ on $\bd U$. By continuity in both variables of the pluricomplex Green function  on a ball there is a ball $W\sbs\sbs U$ such that $dg_B(z,w)>a/2>u(z)-h(z)$ on $\bd
U$ when $w\in W$. For $w\in W$ define the function $v_w(z)$ on $M$ as $g_B(z,w)+h(z)/d$ on $U$, $\max\{g_B(z,w)+h(z)/d,u(z)/d\}$ on $B\sm U$ and $u(z)/d$ outside of $B$. Since
$g_B(z,w)+h(z)/d>u(z)/d$ on $\bd U$ and $g_B(z,w)+h(z)/d<u(z)/d$ on $\bd B$, the function $v_w$ is continuous and plurisubharmonic on $M$. Since $g_M(z,w)\ge v_w(z)=g_B(z,w)+h(z)/d$ on $U$, there is a constant $c$ such that $g_M(z,w)\ge\log\|z-w\|+c$ on $U$ when $w\in W$.
\end{proof}
\subsection{Connections to invariant distances:} The pluri-Greeenian manifolds are Kobayashi hyperbolic, i. e., $k_M(z,w)>-\infty$ for all $z\ne w\in M$ because $k_M(z,w)\ge g_M(z,w)>-\infty$. The converse fails as the example of $\aC\sm\{0,1\}$ shows.
\par Chen and Zhang proved in \cite{CZ} that a manifold with bounded continuous strictly plurisubharmonic functions possesses the Bergman metric, i. e. such manifolds are Bergman hyperbolic.
\par  As the following example shows they need not to be Carath\'eodory hyperbolic, i. e., $c_M(z,w)>-\infty$ for all $z\ne w\in M$, although the class of these manifolds contains all relatively compact domains in Stein manifolds. On such domains the Carath\'eodory function has a strict logarithmic pole at all points.
\bE\label{E:ch} Let $K$ be a regular compact set in $\aC$ such that the Hausdorff measure $H_1(K)=0$ but $H_{1/2}(K)>0$. By \cite[Theorem 5.9.6]{AG} $K$ is not polar. So if $M=\aC\sm K$, then by \cite[Theorem 5.3.8]{AG} $M$ is Greenian. But any bounded holomorphic function on $M$ extends holomorphically to $\aC$ and, consequently, is a constant. Thus $c_M(z,w)\equiv-\infty$.
\eE
\par In \cite{A} Azukawa introduced for the pluri-Greenian manifolds the {\it Azukawa function} $A(w,v)$, where $w\in M$ and $v\in T_M(w)$, as
\[A(w,v)=\limsup_{\zeta\to 0}g_M(f(\zeta),w)-\log|\zeta|,\]
where $f\in\A_r(M)$, $f(0)=w$ and $f'(0)=v$. This definition does not depend on the choice of $f$ but, as Example \ref{E:dcg} shows, $\limsup$ cannot be replaced with $\lim$ and if we put $\liminf$ instead of $\limsup$ then the result will depend on the choice of $f$.
\par In \cite{Ro} Royden introduced the infinitesimal Kobayashi metric that we will use in the form of the {\it Royden function}
\[R_M(w,v)=\inf\{-\log r:\, f\in\A_r(M), f(0)=w, f'(0)=v\}.\]
If $f:\,\aD(0,r)\to M$ is an analytic disk, $f(0)=w$ and $f'(0)=v$, then $g_M(f(\zeta),w)\le g_{\aD(0,r)}(\zeta,0)=\log|\zeta|-\log r$. Thus $A_M(w,v)\le -\log r$. Hence $R_M(w,v)\ge A_M(w,v)$.
\par If $R_M(w,v)=A_M(w,v)$ and there is an extremal mapping $f\in\A_r(M)$ such that $f(0)=w$, $f'(0)=v$ and $-\log r=R_M(w,v)$, then the function $g_M(f(\zeta),w)-\log|\zeta|+\log r$ is subharmonic on $\aD(0,r)$, non-positive and is equal to 0 at 0. Hence $g_M(f(\zeta),w)=\log|\zeta|-\log r$ and we see that it is harmonic on $\aD(0,r)\sm\{0\}$.
\par If $H$ is some Hermitian metric on $M$ and $d_H(z,w)$ is the distance between $z,w\in M$ with respect to this metric, then, instead of comparing $g_M(z,w)$ with logarithms of the Euclidean distance in local coordinates we can compare it with $\log d_H(z,w)$. Clearly a manifold is pluri-Greenian if and only if for every point $w\in M$ there are constants $\si^i_H(w)$ and $\si^s_H(w)$ such that
\[\si_H^i(w)=\liminf_{z\to w} g_M(z,w)-\log d_H(z,w)\le\limsup_{z\to w} g_M(z,w)-\log d_H(z,w)=\si_H^s(w).\]
and if $M$ is locally uniformly pluri-Greenian, then $\si_H^i(w)\ge c$ on a neighborhood of any $w_0\in M$ for some constant $c$ depending on $w_0$.
\par If $f\in\A_r(M)$, $f(0)=w$ and $f'(0)=v$, then
\[\limsup_{\zeta\to 0}g_M(f(\zeta),w))-\log d_H(f(\zeta),w)=A_M(w,v)-\log \|v\|_H=\si_H^s(w,v)\le\si_H^s(w),\]
where $\|v\|_H=H^{1/2}(w,v)$. Since $A_M(w,\al v)=A_M(w,v)+\log|\al|$, $\al\in\aC$, $e^{\si_H^s(w,v)}$ is a well-defined function on $\aCP^{n-1}$ and  can be considered as the {\it directional logarithmic capacity} of $M$ with respect to the metric $H$.
\par On the other hand $A_M(w,v)-\log \|v\|_H\ge\si^i_H(w)$ and we see that
\begin{equation}\label{e:in}e^{\si_H^i(w)}\|v\|_H\le e^{A_M(w,v)}\le e^{\si_H^s(w,v)}\|v\|_H.\end{equation}
\par Since $A_M(w,v)\le R_M(w,v)$ and $R_M(w,v)\le 0$ when there is $f\in\A(M)$ such that $f(0)=w$ and $v=f'(0)$, we obtain the following version of Schwarz Lemma.
\bL\label{T:sl} If $M$ is a pluri-Greenian manifold and $H$ is an Hermitian metric on $M$, then $\|f'(0)\|_H\le  e^{-\si^i_H(w)}$ for any $w\in M$ and any $f\in\A(M)$ with $f(0)=w$. If, additionally, $M$ is locally uniformly pluri-Greenian  and $K$ is a compact set in $M$, then there is a constant $c>0$ such that $\|f'(0)\|_H\le c$ for any $w\in K$ and any $f\in\A(M)$ with $f(0)=w$.
\eL
\par If $B_M(w)$ is the Bergman kernel form on $M$, then we can define a volume form on $M$ that in local coordinates $(z_1,\dots,z_n)$ can be written as $B_M(w)dz_1\wedge d\bar z_1\wedge\cdots\wedge d\bar z_n$. If $dV_H=H(w)dz_1\wedge d\bar z_1\wedge\cdots\wedge d\bar z_n$ is the volume form on $M$ obtained from the metric $H$, then the function $\al_{B,H}(w)=B_M(w)/H(w)$ is well-defined on $M$.
\par If $M$ is a Riemann surface, then $\si^i_H=\si^s_H$. Suita conjectured in \cite{Su} that
\[e^{2\si^s_H(w)}\le\pi\al_{B,H}(w)\]
and the equality holds only if $M$ is the unit disk with a polar set removed. The inequality in Suita's conjecture was proved by B\l ocki in \cite{B1} for bounded domains in $\aC$ and for general Riemann surfaces by Guan and Zhou in \cite{GZ1}. In \cite{GZ2} Guan and Zhou proved the equality part of the conjecture. In \cite{BZ} B\l ocki and Zwonek gave multidimensional generalizations of this conjecture.
\par To get estimates for $\si^i_H$ that appears in Lemma \ref{T:sl} it seems to be reasonable to choose as $H$ the Bergman metric $B(w,v)$. Then the functions $\si^s_B(w)$ and $\si^i_B(w)$ and the form $\si^s_B(w,v)$ are biholomorphically invariant. If $M$ is the unit ball in $\aC^m$ and $0$ is the origin, then $A(0,v)=\log\|v\|_E$, where $E$ is the Euclidean metric. Recall (see \cite[V.19.1 (18) and (19)]{Sh}) that $B(0,v)=(m+1)\|v\|_E^2$, $\|v\|_B=(m+1)^{1/2}\|v\|_E$ and \[d_B(0,z)=(m+1)^{1/2}\log\frac{1+\|z\|_E}{1-\|z\|_E}.\]  Hence $\si^s_B(0)=\si^i_B(0)=-\log(m+1)^{1/2}$ and the inequalities (\ref{e:in}) become identities. The same identities hold at all points of $M$ due to the transitivity of the automorphisms group of $M$.
\par If $M$ is the unit polydisk in $\aC^m$, then using formulas from \cite[V.19.1 (17) and (20)]{Sh} we obtain that $\si^s_B(w)=0$ and $\si^i_B(w)=-\log \sqrt{m}$. It is appealing to conjecture that $\si^i_B\ge-\log(m+1)^{1/2}$ for all complex manifolds with the Bergman metric although the reasons for this are very light.
\subsection{Summary:} As we saw in this section complex manifold with a continuous bounded strictly plurisubharmonic function preserve all properties of relatively compact domains in Stein manifolds, listed in the introduction, except Carath\'eodory hyperbolicity. Boundedness is important: $\aC^n$ has no such function and is not pluri-Greenian.
\par We also would like to indicate that every continuous plurisubharmonic function $u$ on a manifold $M$ with a continuous bounded strictly plurisubharmonic function $v$ can be uniformly approximated by smooth strictly plurisubharmonic functions. Indeed, we may assume that $0\le v\le 1$ and approximate $u$ by functions $u_j=u+2^{-j}v$ and then by Theorem \ref{T:rt} approximate each of $u_j$ by smooth strictly plurisubharmonic functions.
\section{Pluricomplex Green functions on hyperconvex manifolds}\label{S:hcm}
\par The pluricomplex Green functions $g_M(z,w)$ need not to be equal to 0 on the boundaries of bounded domain, for example, when domains are not pseudoconvex and, as the following example shows, they need not to be continuous in $z$ even on bounded pseudoconvex domains $M$ in $\aC^n$.
\bE\label{E:dcg} Take numbers $c_j>0$, converging to 0, and $k_j>0$ such that \[v(z_1,z_2)=\sum\limits_{j=1}^\infty k_j\log|z_2-c_jz_1|\not\equiv-\infty\]
on the ball $\{\|z\|<8\}$, $\sum k_j\log c_j=-\log 2$ and $\sum k_j=1$. Let $u(z)=(\log\|z\|+v(z))/2$ and let $M=\{z\in\aC^2:\,\|z\|<8,\, u(z)<0\}$. Since $g_M(z,0)\ge u(z)$, we see that $g_M((1/2,0),0)\ge-\log 2^{3/2}$. On the other hand the points $z_j=(1/2,c_j/2)$ lie on the intersections of the lines $z_2=c_jz_1$ with $M$ that are disks of radius 8 centered at 0 and converge to $(1/2,0)$. Thus \[\limsup_{j\to\infty}g_M(z_j,0)\le\limsup_{j\to\infty}\frac12\log((1+c_j^2)^{1/2}/16)
=-2\log2<-\log 2^{3/2}\]
and this shows the discontinuity of $g_M$ in $z$.
\par Also this example shows that one cannot define an analog of the Azukawa function via $\liminf$. Indeed, if $f(\zeta)=(\zeta,0)\in\A_r(M)$, then
\[\liminf_{\zeta\to0}g_M(f(\zeta),0)-\log|\zeta|\ge
\liminf_{\zeta\to0}u(f(\zeta))-\log|\zeta|=-\frac12\log 2\]
If $g(\zeta)=(\zeta,\zeta^2)\in\A_r(M)$, $z_j=(c_j,c_j^2)$ and $\zeta_j=c_j$, then
\[\liminf_{\zeta\to0}g_M(g(\zeta),0)-\log|\zeta|
\le\liminf_{j\to\infty}g_{\aD(0,8)}(\zeta_j,0)-\log|\zeta_j|\le-3\log 2\]
and we see that $\liminf$ depends on the choice of $f$.\eE
\par For bounded domains both properties  of $g_M(z,w)$ mentioned above are guaranteed by hyperconvexity. A bounded domain $M$ is called {\it hyperconvex} if there is a negative plurisubharmonic function $u$ on $M$ such that for any $r<0$ the set $B_{u,r}=\{z\in M:\,u(z)<r\}$ is relatively compact in $M$. We will call $u$ the {\it negative exhaustion function}. It was proved by Kerzman and Rosay in \cite{KR} that on hyperconvex domains in $\aC^n$ there are smooth strictly plurisubharmonic negative exhaustion functions. In \cite[Theorem 4.3]{D1} Demailly has shown that the pluricomplex Green function $g_M(z,w)$ on a relatively compact hyperconvex domain in a Stein manifold is continuous on $\ovr M\times M$.
\par In \cite{C1} Chen introduced hyperconvex manifolds as manifolds with smooth negative strictly plurisubharmonic exhaustion function and proved that the Bergman metric is complete on such manifolds. We relax this definition a bit but show later that our definition is equivalent to Chen's definition.
\par We say that a complex manifold $M$ is {\it hyperconvex} if it has a negative continuous plurisubharmonic exhaustion function and a bounded continuous strictly plurisubharmonic function. By a theorem of Narasimhan (see \cite[Theorem II]{Nar}) a manifold with a continuous plurisubharmonic exhaustion function and a strictly plurisubharmonic function is Stein.
\par The need for the second function can be explained by the following example from \cite{PSh}.
\bE Take the unit ball $B$ in $\aC^2$ and blow-up a complex projective line $X$ at the origin. We get a complex manifold $M$ and a holomorphic mapping $F$ of $M$ onto $B$ such that $F(X)=\{0\}$. If $u$ is a plurisubharmonic function on
$B$, then the function $v=u\circ F$ is plurisubharmonic on $M$ while
any plurisubharmonic function on $M$ is constant on $X$. Hence $M$ has a negative exhaustion function but for any $w\in X$ the function $g_M(z,w)\equiv-\infty$ on $X$.
\eE
\par Now we will show that on hyperconvex manifolds the pluricomplex Green functions have the properties listed in Section \ref{S:1}.
\bT\label{T:pcgh} Let $M$ be a hyperconvex manifold. Then:
\be\item $M$ is locally uniformly pluri-Greenian;
\item for any $w\in M$ the function $g_M(z,w)$ is a negative exhaustion function;
\item the function $g_M(z,w)$ is continuous  on $\ovr M\times M$;
\item $M$ has a smooth strictly plurisubharmonic negative exhaustion function.
\ee
\eT
\par{\bf Remark:} Item (4) implies that our definition is equivalent to Chen's definition of hyperconvex manifolds in \cite{C1}.
\begin{proof} (1) follows immediately from Theorem \ref{T:bsplug} and this implies that for any $w\in M$ the function $g_M(z,w)$ is locally bounded on $M\sm\{w\}$ and, by its maximality, the sets $M_a=\{z\in M:\,g_M(z,w)<a\}$ are connected for any $a<0$.
\par To show (2) we will follow the proof of Theorem 4.3 in \cite{D1}. Let $u$ be a negative plurisubharmonic exhaustion function on $M$. Let us take a coordinate neighborhood $U$ of $w\in M$. Since $M$ is pluri-Greenian, there is an open ball $B$ of radius $r$ centered at $w$ such that $g_M(z,w)>\log r+c$ on $\bd B$. Let $\al$ be the maximum of $u$ on $\bd B$. Clearly $\al<0$ and there is a constant $b>0$ such that $bu<g_M(z,w)$ on $\bd B$. Hence $g_M(z,w)\ge bu(z)$ on $M\sm B$ by the maximality of $g_M$ and (2) is proved.
\par Let $\{(z_j,w_j)\}$ be a sequence  converging to $(z_0,w_0)$ in $M\times M$. If $z_0=w_0$, then $g_M(z_j,w_j)\to -\infty=g_M(w_0,w_0)$ because  $M$ is locally uniformly pluri-Greenian. So we may assume that $z_0\ne w_0$. By Lemma \ref{L:kl} there is a sequence of $\eps_j>0$ converging to 0 such that
\[(1+\eps_j)g_M(z_j,w_0)\le g_M(z_j,w_j)\le (1-\eps_j)g_M(z_j,w_0).\]
Therefore to prove that the sequence $\{g_M(z_j,w_j)\}$ converges to $g_M(z_0,w_0)$ it suffices to prove that the sequence $\{g_M(z_j,w_0)\}$ converges to $g_M(z_0,w_0)$. Since $g_M(z,w_0)$ is upper semicontinuous $\limsup_{j\to\infty}g_M(z_j,w_0)\le g_M(z_0,w_0)$.
\par To prove that $\liminf_{j\to\infty}g_M(z_j,w_0)\ge g_M(z_0,w_0)$ we cannot follow \cite{D1} because that proof uses the result of Walsh in \cite{W} that requires an attainable boundary and it is missing in the case of manifolds. We will use the definition of the pluricomplex Green functions via analytic disks.
\par Let us take a sequence $\{\eps_j\}$ of negative real numbers strictly increasing to 0 and let $M_j=\{z\in M:\,g_M(z,w_0)<\eps_j\}$. The numbers $\dl_j=\liminf_{z\to\bd M_j}g_M(z,w_0)$ also converge to 0 and the open set $M'_j=\{z\in M:\,g_M(z,w_0)<2\dl_j\}$ is relatively compact in $M_j$.
\par Clearly, $g_M(z,w_0)\le g_{M_j}(z,w_0)$ on $M_j$. On the other hand, let us take any $\al>0$. Since all manifolds above are pluri-Greenian, for each $j$ there is a neighborhood $U_j$ of $w_0$ such that $(1+\al)g_{M_j}(z,w_0)<g_M(z,w_0)-2\dl_j$ on $U_j$. Moreover, $(1+\al)g_{M_j}(z,w_0)<g_M(z,w_0)-2\dl_j$ on $\bd M'_j$ because $g_M(z,w_0)-2\dl_j\ge 0$ there. Hence $(1+\al)g_{M_j}(z,w_0)<g_M(z,w_0)-2\dl_j$ on $M_j'$ by the maximality of $g_M$. Since $\al$ can be chosen arbitrarily small we see that  $g_{M_j}(z,w_0)\le g_M(z,w_0)-2\dl_j$ on $M_j'\sm\{w_0\}$. It follows that the functions $g_{M_j}(z,w_0)$ converge to $g_M(z,w_0)$ uniformly on compacta in $M\sm\{w_0\}$.
\par Since $M$ is Stein, there is an imbedding $F$ of $M$ into $\aC^N$. By \cite[Theorem 8.C.8]{GR} there are an open neighborhood $U_j$ of $F(\ovr M_j)$ in ${\mathbb C}^N$ and a holomorphic retraction $P_j$ of $U_j$ onto $F(M)$.
\par Let $z'=F(z)$. There are invertible affine transformations $A_k$ of $\aC^N$ such that $A_k(z'_k)=z'_0$, $A_k(w'_0)=w'_0$ and the transformations $A_k$  converge to the identity mapping uniformly on compacta in $\aC^N$. Let us choose some $\eps>0$ and find an analytic disk $f_{jk}:\,\oD\to M_j$ such that $f_{jk}(0)=z_k$ and $H_{w_0}(f_{jk})<g_{M_j}(z_k,w_0)+\eps$. For sufficiently large $k$ the disks $\wtl g_{jk}=A_k(F(f_{jk}))$ lie in $U_j$, $\wtl
g_{jk}(0)=z'_0$ and $\wtl g_{jk}(\zeta)=w'_0$ when $f_{jk}(\zeta)=w_0$. Hence for the disks $g_{jk}=P_j(\wtl g_{jk})$ the functional $H_{w_0}(g_{jk})=H_{w_0}(f_{jk})<g_{M_j}(z_k,w_0)+\eps$ and, therefore,  $g_M(z_0,w_0)<g_{M_j}(z_k,w_0)+\eps$.
\par Since the functions $g_{M_j}(z,w_0)$ converge to $g_M(z,w_0)$ uniformly on compacta in $M\sm\{w_0\}$, for any $\dl>0$ there is $j$ such that  $g_{M_j}(z_k,w_0)<g_{M}(z_k,w_0)+\dl$ for all $k$. Thus
\[g_M(z_0,w_0)\le\liminf_{k\to\infty}g_{M_j}(z_k,w_0)+\eps\le\liminf_{k\to\infty} g_{M}(z_k,w_0)+\eps+\dl.\]
Taking into account that $\eps>0$ and $\dl>0$ are arbitrary we see that
\[\liminf_{k\to\infty}g_M(z_k,w_0)\ge g_M(z_0,w_0).\]
Hence $\lim_{k\to\infty}g_M(z_k,w_k)=g_M(z_0,w_0)$ and (3) is proved.
\par Let $u$ be a bounded continuous strictly plurisubharmonic function on $M$. We may assume that $-1\le u<0$ on $M$. Take a point $w_0\in M$ and define $u_j(z)=\max\{jg_M(z,w_0),u(z)\}$. The function $v(z)=\sum_{j=1}^\infty2^{-j}u_j(z)$ is a continuous strictly plurisubharmonic exhaustion function on $M$. In Theorem \ref{T:rt} we take the function $\lm=-v/2$ and find a smooth strictly plurisubharmonic function $\wtl v$ on $M$ such that $v\le\wtl u\le v/2$. Hence $u$ is smooth strictly plurisubharmonic exhaustion function on $M$ and this proves (4).
\end{proof}
\subsection{Connections to invariant distances:}
\par As Example \ref{E:ch} shows hyperconvex manifolds need not to be Carath\'eodory hyperbolic because a domain in $\aC$ with the regular non-polar complement is hyperconvex. However, hyperconvex manifolds are not only Kobayashi hyperbolic but also Kobayashi complete because $k_M(z,w)\ge g_M(z,w)$ and all Kobayashi balls compactly belong to $M$. In \cite{C1} Chen proved that hyperconvex manifolds are Bergman complete.
\par Arguing as in the proof of the continuity of the pluricomplex Green functions, one can show that the Kobayashi and Royden functions are continuous on $M\times M$ and $T_M$ respectively. Note that if $f\in\A_r(M)$, then $g_M(f(\zeta),f(0))\le\log(|\zeta|/r)$. From this it is straightforward to prove that the Kobayashi and Royden variational problems attain their extrema.
\subsection{Summary:}  As we saw in this section hyperconvex complex manifold preserve all basic  properties of hyperconvex relatively compact domains in Stein manifolds, listed in the introduction.

\end{document}